\documentclass{aes2e}
\usepackage{soul}
\usepackage{float}
\usepackage{algorithm}
\usepackage[noend]{algpseudocode}
\usepackage[T1]{fontenc}
\usepackage{newtxmath}
\usepackage{euscript}
\usepackage{graphicx}
\usepackage{mathtools}

\DeclarePairedDelimiter\floor{\lfloor}{\rfloor}
\graphicspath{ {./images/} }

\usepackage{hyperref}
\newcommand\blfootnote[1]{%
  \begingroup
  \renewcommand\thefootnote{}\footnote{#1}%
  \addtocounter{footnote}{-1}%
  \endgroup
}

\makeatletter
\def\algbackskip{\hskip-\ALG@thistlm}
\makeatother

\jyear{2020}
\jmonth{November}

\begin{document}


\title{\emph{GL-Coarsener:} A Graph representation learning framework to construct coarse grid hierarchy for AMG solvers}

\authorgroup{
\author{Reza Namazi}\footnotemark,\quad
\author{Arsham Zolanvari}\footnotemark,\quad
\author{Mahdi Sani}\footnotemark,\quad
\author{Seyed Amir Ali Ghafourian Ghahramani}\footnotemark
}
\abstract{%
 In many numerical schemes, the computational complexity scales non-linearly with the problem size. Solving a linear system of equations using direct methods or most iterative methods is a typical example. Algebraic multi-grid (AMG) methods are numerical methods used to solve large linear systems of equations efficiently. One of the main differences between AMG methods is how the coarser grid is constructed from a given fine grid. There are two main classes of AMG methods; graph and aggregation based coarsening methods. Here we propose an aggregation-based coarsening framework leveraging graph representation learning and clustering algorithms. Our method introduces the power of machine learning into the AMG research field and opens a new perspective for future researches. The proposed method uses graph representation learning techniques to learn latent features of the graph obtained from the underlying matrix of coefficients. Using these extracted features, we generated a coarser grid from the fine grid. The proposed method is highly capable of parallel computations. Our experiments show that the proposed method's efficiency in solving large systems is closely comparable with other aggregation-based methods, demonstrating the high capability of graph representation learning in designing multi-grid solvers.
}
\maketitle
\blfootnote{*\hspace{1pt}\textbf{Reza Namazi} is a Computer Engineering B.Sc. Student at Sharif University of Technology - Kish International Campus, Email: \href{mailto:rezanmz@ymail.com}{rezanmz@ymail.com}}
\blfootnote{†\hspace{1pt}\textbf{Arsham Zolanvari} is an Industrial Engineering B.Sc. Student at Sharif University of Technology - Kish International Campus, Email: \href{mailto:arsham.zollanvari@gmail.com}{arsham.zollanvari@gmail.com}}
\blfootnote{‡\hspace{1pt}\textbf{Mahdi Sani} is an Assistant Professor of Mechanical Engineering at Sharif University of Technology - Kish International Campus, Email: \href{mailto:msani@sharif.edu}{msani@sharif.edu}}
\blfootnote{§\hspace{1pt}\textbf{Seyed Amir Ali Ghafourian Ghahramani} is an Assistant Professor of Computer Engineering at Sharif University of Technology - Kish International Campus, \quad\quad\quad\quad\quad\quad\quad\quad\quad Email: \href{mailto:ghahramani@ce.sharif.edu}{ghahramani@ce.sharif.edu}}
\section{INTRODUCTION}
Many real-world problems are governed by
\emph{partial differential equations (PDEs)}. For example, to predict the forces induced by airflow on a car or an airplane, Computational Fluid Dynamics (CFD) relies on the solution of non-linear Navier Stokes Equations (NSE). To model these problems, large complex meshes are constructed. After applying numerical methods, it results in a massive system of linear equations (order of million to billion) to be solved many times (order of thousand to million) during the simulation process \cite{stuben2001review}.
\par
The typical way to solve PDEs is to discretize the PDE to equations that involve a finite number of unknowns. This is usually achieved using \emph{Finite Differences Method (FDM)}, \emph{Finite Volume Method (FVM)}, or \emph{Finite Element Method (FEM)}. At the heart of these processes is solving efficiently a linear systems of equations in the form of:
\begin{equation}
    Au = f
    \label{eq:system}
\end{equation}
where $A$ is a sparse coefficient matrix resulted from the discretization of the original PDE which is called \emph{stiffness matrix} in some cases.
\par
There are many methods that could be used to solve systems of linear equations similar to \textbf{Eq. \ref{eq:system}}. These methods could be sorted into two main categories;
\emph{direct} or \emph{iterative} methods.
\par
\emph{Direct} methods like \emph{variable elimination techniques}, \emph{row reduction techniques (e.g., Gaussian elimination)}, \emph{Cramer's technique} and \emph{inverse matrix
solution} will result in an exact (up to machine accuracy) solution. Direct methods have a typical time complexity of $\mathcal{O}(n^3)$ \cite{davis2016survey}. As the system's size increases, it becomes computationally costly to use \emph{direct} methods, which is the main reason they are not usually used for real-world engineering problems. Moreover, they usually destroy the original matrix's sparsity, making it hard to keep on high-speed \emph{Random Access Memory (RAM)}.
\par
\emph{Indirect (iterative)} methods like \emph{Jacobi}, \emph{Gauss-Seidel}, and \emph{Successive Over-Relaxation (SOR)} usually converge to an accurate enough approximation for the system \cite{cerwinsky2013theory}. Under the right condition, these iterative methods will rapidly damp the approximate solution's high-frequency errors in a few steps -\emph{smoothing the approximate solution}- but are ineffective when it comes to damping low-frequency errors. Two important properties of the iterative methods are that they provide good approximations in manageable amounts of time and they do not destroy the sparsity of the coefficient matrix.
\par
\emph{Multi-grid (MG)} methods are a powerful class of iterative solvers that could efficiently solve discretized differential equations. As mentioned previously, \emph{iterative solvers (smoothers)} like \emph{Jacobi} are useful for damping high-frequency errors. Unfortunately, after a few iterations, their convergence rate decreases drasically and usually the slow convergence rate associated with the low-frequency errors dominates the error reduction rate. To address this issue, \emph{multi-grid} methods introduce a hierarchy of coarser grids. With proper mapping (called prolongation) the low-frequency errors with respect to the fine grid are represented as high-frequency error relative to coarse grid; therefore the smoothers would keep their high rate of convergence on the coarse grid. As such, multi-grid methods can be regarded more as a strategy than stand-alone solvers \cite{kepleralgebraic}.
\par
A coarse grid is a lower resolution representation of the original fine grid. The key to success of multi-grid strategy is the fact that the low-frequency errors of the finer grid will become high frequency in the coarser grid \cite{strang2006mathematical}.
\par
The coarse grid problem could either be constructed from the original physics (Physical or Geometric Multi Grid) or constructed from the fine grid's existing coefficient matrix (Algebraic Multi-grid, \emph{AMG}). In this work, we focus on \emph{AMG}. One of the main differences between different \emph{AMG} methods is how they move from a fine grid to a coarser grid. When it comes to constructing the coarser grid from the fine grid, \emph{AMG} methods fall into two main categories \cite{garcia2019parallel}:
\par
\emph{1-} Graph based AMG methods \emph{(e.g., Beck, Ruge-Stüben)}
\par
\emph{2-} Aggregation based AMG methods \emph{(e.g., Vaněk, Adaptive Smoothed
Aggregation ($\alpha SA)$)}
\par
This paper proposes \emph{\textbf{GL-Coarsener}} framework, an aggregation-based \emph{AMG}, for constructing the coarser grid linear system from that of the given fine grid leveraging the underlying connections in the fine grid's graph.
\newline
\emph{\textbf{GL-Coarsener}} consists of two main modules:
\par
\emph{- Embedding module}: The underlying fine grid is projected into a \emph{d-dimensional} embedding space.
\par
\emph{- Clustering module}: The aggregates are constructed by applying clustering techniques on the embedding space obtained from the \emph{embedding module}.
\par
To the best of our knowledge, our work is the first \emph{aggregation based AMG}
that utilizes graph representation learning techniques to select aggregates intelligently.
\par
The rest of the paper is organized as follows, In section 2 we briefly review the related works and researches done in \emph{AMG} field. In section 3 the required preliminaries will be explained. We present our proposed method in section 4. Finally, in section 5 we compare the performance of our method with other available methods.
\section{RELATED WORK}
Algebraic multi-grid was first described  in 1980s. \cite{brandt1985algebraic}. This method assumes no information is provided about the geometry of the underlying grid and takes the \textbf{coefficient matrix}  $(A)$ which could be obtained from the discretized PDE problem to construct the multilevel hierarchy. The earliest algebraic multi-grid methods are graph-based methods. By introducing some criteria to determine strongly connected nodes, graph-based methods label nodes as \emph{coarse} or \emph{fine} by taking into account the count of a node's strong connections and the strength of connections. 
\par
A popular graph-based algorithm is Ruge-Stüben \cite{ruge1987algebraic}, introduced in 1981, this method and it's refined versions are still widely used. \emph{Ruge-Stüben} takes a threshold to distinguish between weak and strong connections; the nodes with more strong connections are selected to be coarse nodes. The prolongation operator is then constructed by considering the number of strong connections of each node. Later implementations of this method include P. Zaspel \cite{zaspel2017analysis}, which proposes a parallelized GPU implementation of \emph{Ruge-Stüben} to improve the algorithm's performance.
\par
The \emph{Beck} method is a simpler coarsening method which does not distinguish between strong and weak connections \cite{beck1999graph}. This results in less robustness, meaning that it cannot recover the usual multi-grid efficiency for non-smooth coefficient PDEs, but allows for an easier and less complex implementation. \emph{Beck} method carries out the task of labelling nodes and constructs transfer operators by considering only the number of connections.
\par
Aggregation\footnote{\emph{Aggregation} is used in two different contexts in this article. One refers to  grouping  nodes in AMG context and the other refers to accumulation of node features in graph learning algorithms.} based methods are a class of coarseners that differ from graph based methods in how the coarse grid is constructed. Aggregate methods make use of connection strength, but only to form neighborhoods of points. These neighborhoods are then taken as a unit while selecting the coarse grid \cite{cerwinsky2013theory}. The first aggregation method was introduced by Vaněk. The standard aggregation \cite{vanvek1996algebraic} defines the strongly coupled neighborhood of a node \textbf{i} with  threshold {$\epsilon$}  as follows:
\begin{equation}
    N_{i }(\epsilon ) = \left \{ j\in \Lambda _{h} : \left | A_{ij}  \right  | \geq \epsilon  \sqrt{\left | A_{ii}A_{jj}^{} \right |} \right \}
\end{equation}
where $\Lambda_h$ is the set of node indices in fine level adjacency matrix.
These neighborhoods are used to form aggregates of nodes denoted by \textbf{C}. Once the aggregates are selected, the corresponding prolongation matrix is constructed using the cluster (aggregate) assigned to each node.
\par
The preliminary prolongation matrix \textbf{$\hat{P}$} in standard aggregation is defined as the following:
\begin{equation}
    \mathbf{\hat{P}}_{ij} = \left\{\begin{matrix}
 1,\: \: \: \: \: \: \: \: \:  i \in \mathbf{C_{\mathit{j}}};
\\
0, \: \: otherwise.
\end{matrix}\right.
\label{eq:prolongator}
\end{equation}
where $i$ and $j$ iterate over nodes and clusters respectively. This prolongation operator is used in this paper.
\par
Recent improvements on aggregation based methods include adaptive aggregation \cite{de2008multilevel}, and accelerated adaptive aggregation methods \cite{shen2018block}.
\begin{figure*}[ht]
    \centering
    \includegraphics[height=8.55cm]{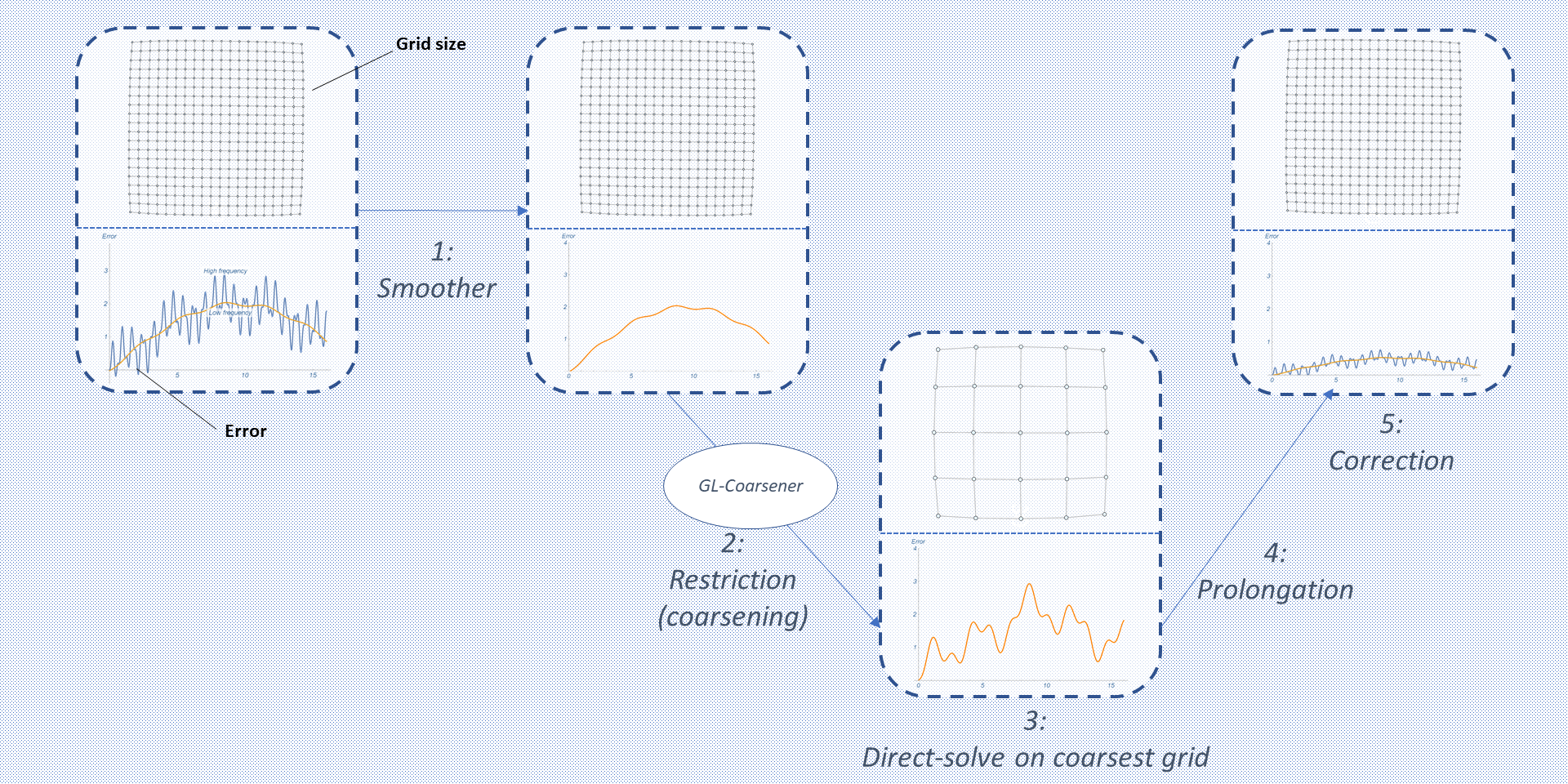}
    \caption{Illustration of a two-grid V-cycle}
    \label{fig:multigrid_flow}
\end{figure*}
\section{PRELIMINARIES}
In this section, we first review the basics of \emph{AMG}, then we describe the machine learning techniques used in this paper.
\subsection{Algebraic Multi-grid (AMG)}
It is well known that stationary iterative linear solvers \emph{(like Jacobi and Gauss-Seidel)} cannot efficiently damp low-frequency errors, especially for large sparse symmetric positive definite (SPD) systems. Further analysis of iterative solvers has proven that this loss of convergence rate is related to the system's smallest eigenvalues. To recover the convergence rate, algebraic multi-grid methods construct a hierarchy of operators to map the system of equations to a coarser system with a fewer number of unknowns. As is shown in \textbf{Fig. \ref{fig:multigrid_flow}}, the low-frequency errors of the fine grid problem (related to lowest eigenvalues) are represented as high-frequency errors on the next coarse grid. Progressively constructing coarser grids, \emph{AMG} methods will construct systems of equations with relative high eigenvalues; therefore, the iterative smoothers' convergence rate is recovered on coarse grid and it could efficiently reduce the errors again.
\\\\
As demonstrated in \textbf{Fig. \ref{fig:multigrid_flow}}, \emph{AMG} methods work in five main phases:
\par
\emph{1-} The approximate solution $v$ obtained from applying a few iterations of the \emph{smoother} will be used to calculate the fine grid residual:
\begin{equation}
    \begin{gathered}
        r = f - Av
        \\
        = Au - Av
        \\
        = A(u - v) = Ae
    \end{gathered}
    \label{eq:residual}
\end{equation}
in which $f$ is the right-hand-side matrix of the original system (\textbf{Eq. \ref{eq:system}}) and $e$ is the error involved in approximating $u$ with $v$.
\par
\emph{2-} A coarse grid is selected, and a restriction operator is constructed to map the \emph{coefficient matrix ($A_f$)} from the fine grid to the coarse grid.
\par
\emph{3-} The residual \textbf{Eq. \ref{eq:residual}} is solved on the coarse grid. The solution in this step is accelerated by recursive application of the multi-grid method until the coarsest grid is reached, for which obtaining an exact solution using direct methods would be cheap.
\par
\emph{4-} The third phase's solution is interpolated back to the fine grid, usually using the restriction operator's transpose.
\par
\emph{5-} The interpolated solution is used to modify the approximate solution $v$ to  $v+e$. This correction will reduce low-frequency errors in $v$ for which the smoother had a low rate of convergence. Before and after each multi-grid iteration, a smoothing step is usually applied to the system of equations to damp the high-frequency errors. These steps are called \emph{pre-smoothing} and \emph{post-smoothing}. \footnote{Usually a small number of iterations (1 to 5) of smoothing is sufficient}
\\
AMG algorithms are usually described as a recursive application of the two-grid method. \textbf{Algorithm \ref{alg:general_multigrid}} explains a two-level AMG algorithm.
\begin{algorithm}[H]
\caption{One iter. of two-level Algebraic Multi-grid}
\hspace*{\algorithmicindent} \textbf{Input:} $A_f, f_f$ \algorithmiccomment{$A_fu_f=f_f$ is the fine-level system}\\
\hspace*{\algorithmicindent} \hspace{1.1cm}$v_{f}$ \algorithmiccomment{$v_f$ is an approximation of $u_f$}\\
\hspace*{\algorithmicindent} \textbf{Output:} \emph{Corrected $v_{f}$}
\begin{algorithmic}[1]
\State $\emph{$r_{f}$} \gets \emph{$f_{f}$ - $A_{f}v_{f}$}$
\State Construct$\emph{ prolongation operator \textbf{P}}$
\State $\emph{R} \gets \emph{P}^T$
\State $A_{c} \gets RAP$
\State $r_{c} \gets Rr_{f}$
\State Solve $A_{c}e_{c} = r_{c}$ recursively for $e_{c}$
\State $e_{f} \gets Pr_{c}$
\State \textbf{return} $v_{f} + e_{f}$
\end{algorithmic}
\label{alg:general_multigrid}
\end{algorithm}
\begin{figure*}[ht]
    \centering
    \includegraphics[height=10.3cm]{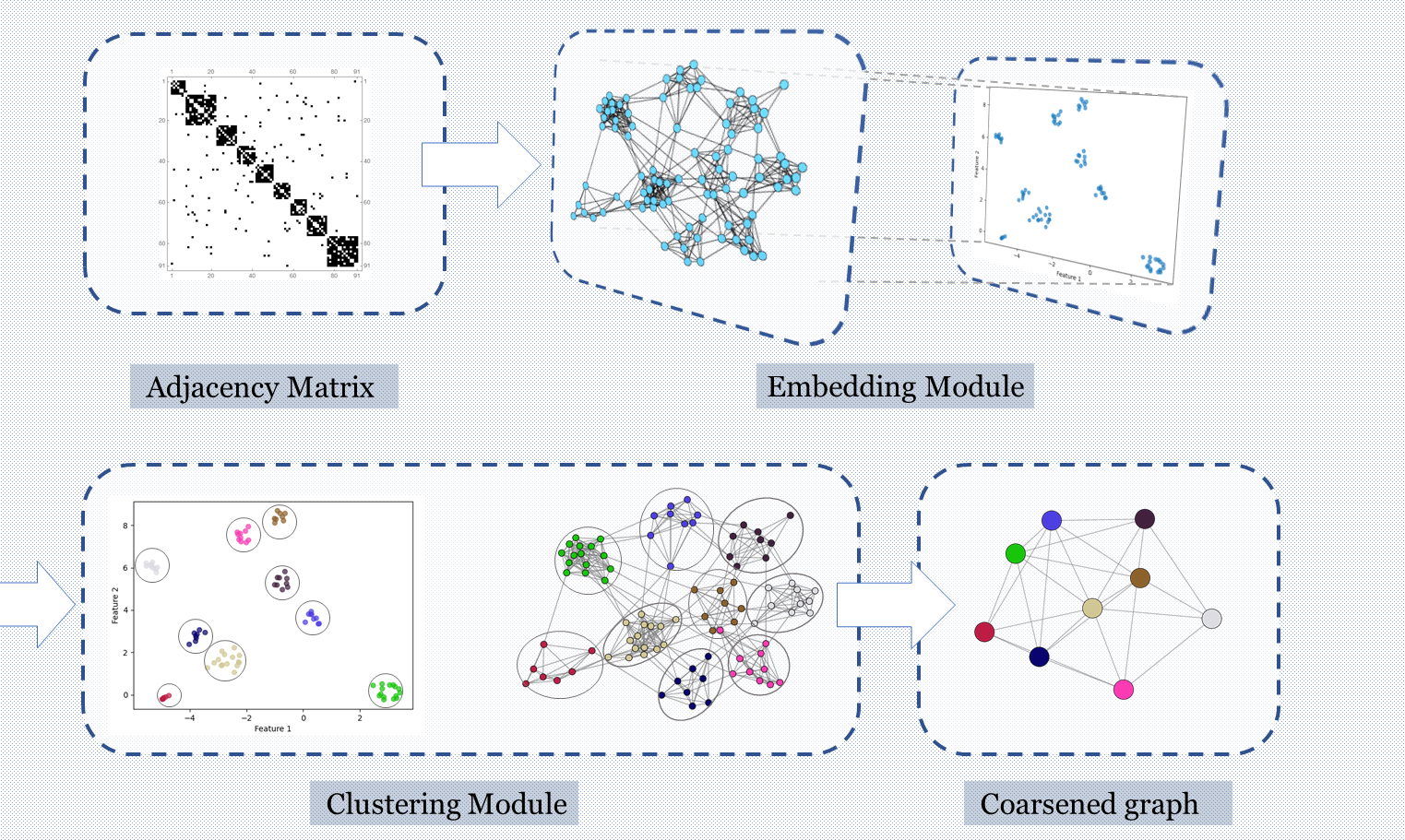}
    \caption{$GL-Coarsener$}
    \label{fig:GL_coarsener}
\end{figure*}
\subsection{Machine Learning}
\noindent
\subsubsection{node2vec}
Graph learning and feature representation are the main concepts behind the \emph{embedding module} of our method. This would allow for accurate downstream machine learning tasks such as clustering.
\par
Dimensionality reduction techniques for unsupervised feature learning such as Vazquez \cite{vazquez2003global} and Tenenbaum \cite{tenenbaum2000global} could be used in this work, but these methods show lower computational performance, and become increasingly more difficult when scaled to large graphs. The Skip-gram model \cite{mikolov2013distributed} uses neural network models to learn word feature
vectors from huge data sets with billions of words. This algorithm is  based on the hypothesis that similar words tend to be closer to each other. Moreover, words can have multiple degrees of similarity \cite{mikolov2013linguistic}. The \emph{skip-gram} model and the idea of feature representation was further extended to networks by \emph{DeepWalk} algorithm \cite{perozzi2014deepwalk}; \emph{DeepWalk} generates random-walks on nodes to analyze graphs and learn latent representations of vertices in a network . This work is further developed by \emph{node2vec} algorithm. The algorithm adopts a flexible parameter tuning for node sampling  to learn richer representation in different networks \cite{grover2016node2vec}.
\par
One of the newer methods for representing graphs is using convolutional neural networks (CNN) \cite{niepert2016learning}. The basic idea is to slide a filter over the structure of the grid to learn representation of nodes. However, this approach makes the algorithm dependent on the structure of the grid therefore the training does not generalize to various formations of graphs. Graph convolutional networks (GCN), on the other hand, learn aggregations of a node's neighbours to construct the node's feature vector \cite{kipf2016semi}. Such methods are independent of the ordering of nodes and the structure of the graph. This grants the embedder \emph{transductivity} property, which allows for generation of node embeddings to previously unseen data. GraphSAGE algorithm \cite{hamilton2017inductive} is a developed GCN variant, which trains a set of \emph{aggregator functions} that learn to aggregate feature information from a node’s local neighborhood. GraphSAGE utilizes three aggregator functions; Mean, LSTM and pooling aggregators.
\par
\emph{\textbf{GL-Coarsener}} represents the \emph{coefficient matrix} as a graph, enabling the implementation of state-of-the-art embedding algorithms such as those mentioned in this section. A precise graph embedding plays an essential role in an effective coarsening.

Learning graph latent features has always been an interesting topic in network analysis field. In this paper, we use \emph{node2vec} \cite{grover2016node2vec}, a graph representation learning algorithm that learns node features in a graph using machine learning techniques. 
\par
The goal of graph representation learning is to map each node of the network \emph{G(V, E)} to a \emph{d-dimensional} vector in the embedding space $Z$ in a way that similar nodes in the graph are close to each other in the \emph{d-dimensional} embedding space. This allows us to encode network information and generate node representation. The \emph{Embedding Module} in \textbf{Fig. \ref{fig:GL_coarsener}} illustrates a two-dimensional representation of a graph.
\par
We can use different metrics to measure the similarity of the nodes in the embedding space; here we are using cosine similarity (dot product) \cite{grover2016node2vec}:
\begin{equation}
    similarity(u, v) \approx z_v^Tz_u
\end{equation}
For measuring the similarity of the nodes in the original graph we will use \emph{random walk} approach. We will explore the graph by applying $\gamma$ random walks of length $t$ on each node and define the similarity of nodes \emph{u} and \emph{v} as the probability that \emph{u} and \emph{v} co-occur on a random walk over the network.
\par
Each random walk generates a set of neighbors for each node. For example, for node $u$ random walk approach generates $\EuScript{N}(u)$.
\par
To generate the embedding vectors, we can view this as an optimization problem with the goal of maximizing the similarity of the neighboring nodes obtained from random walks \cite{grover2016node2vec}:
\begin{equation}
    max\prod_{u \in V}P(\EuScript{N}(u)|z_u)
    \label{eq:embedding_probability}
\end{equation}
To obtain the embedding space similarity, \textbf{Eq. \ref{eq:embedding_probability}} calculates the product of all node similarities. Multiplying the probabilities, will cause the result to vanish quickly. To address this phenomenon, we take the $\log$ of \textbf{Eq. \ref{eq:embedding_probability}}:
\begin{equation}
    \begin{gathered}
        max(\log\prod_{u \in V}P(\EuScript{N}(u)|z_u))
        \\
        = max\sum_{u \in V}\log P(\EuScript{N}(u)|z_u)
    \end{gathered}
\end{equation}
The optimization problem could be redefined as minimizing the loss function:
\begin{equation}
    \mathcal{L} = \sum_{u \in V} \sum_{v \in \EuScript{N}(u)} - \log{P(v|z_u)}
\end{equation}
where $P(v|z_u)$ is the probability of reaching node $v$ in the random walks given that we started the walk from node $u$. By minimizing the above equation, we get an embedding space that nodes $u$ and $v$ are closer if they are neighbors in the random walks obtained from the graph.
\par
To parameterize $P(v|z_u)$, we will use softmax \cite{gao2017properties}:
\begin{equation}
    P(v|z_u) = \frac{exp(z_u^Tz_v)}{\sum_{n \in V} exp(z_u^Tz_n)}
\end{equation}
Softmax ensures that the output is between 0 and 1. The equation above is simply used so that node $v$ be most similar to node $u$ out of all nodes $n$.
\par
Putting it all together, we have:
\begin{equation}
    \mathcal{L} = \sum_{u \in V} \sum_{v \in \EuScript{N}(u)} - \log{\frac{exp(z_u^Tz_v)}{\sum_{n \in V} exp(z_u^Tz_n)}}
\end{equation}
Optimizing the random walks basically means finding the embeddings $z_u$ that minimizes $\mathcal{L}$.
\par
The problem with the above equation is that nested sum over all nodes in the graph is computationally very expensive, it gives a computational complexity of $\mathcal{O}(|V|^2)$. To address this issue, the sum over all nodes in the softmax equation could be approximated as \cite{grover2016node2vec}:
\begin{equation}
    \log{\frac{exp(z_u^Tz_v)}{\sum_{n \in V} exp(z_u^Tz_n)}} \approx
    \log{\sigma(z_u^Tz_v)} - \sum_{i=1}^{k}\log{\sigma(z_u^Tz_{n_i})}
\end{equation}
In the above equation, $\sigma$ is the sigmoid function and $n_i$ is selected randomly from all nodes in the graph. This technique is called \emph{\textbf{negative sampling}} \cite{goldberg2014word2vec}. Instead of normalizing the cosine similarity with respect to all nodes, we just normalize against $k$ random negative samples $n_i$. Higher $k$ gives a more robust estimates.
\par
Now that the general optimization based algorithm was explained, we will discuss two different methods that could be used for generating random walks; \emph{\textbf{DeepWalk}} and \emph{\textbf{node2vec}}
\par
\emph{\textbf{DeepWalk}} is the simplest random walk approach. For generating the random walks, it will start from a given node and selects the next node randomly based on uniform distribution, generating an unbiased random walk.
\par
The next method is \emph{\textbf{node2vec}}. The idea of \emph{\textbf{node2vec}} is to generate a biased random walk that can trade off between local (micro) and global (macro) views of the graph. It is done by biased $2^{nd}$-order random walks that explore network neighborhoods.
\begin{figure}[H]
    \centering
    \includegraphics[height=5cm]{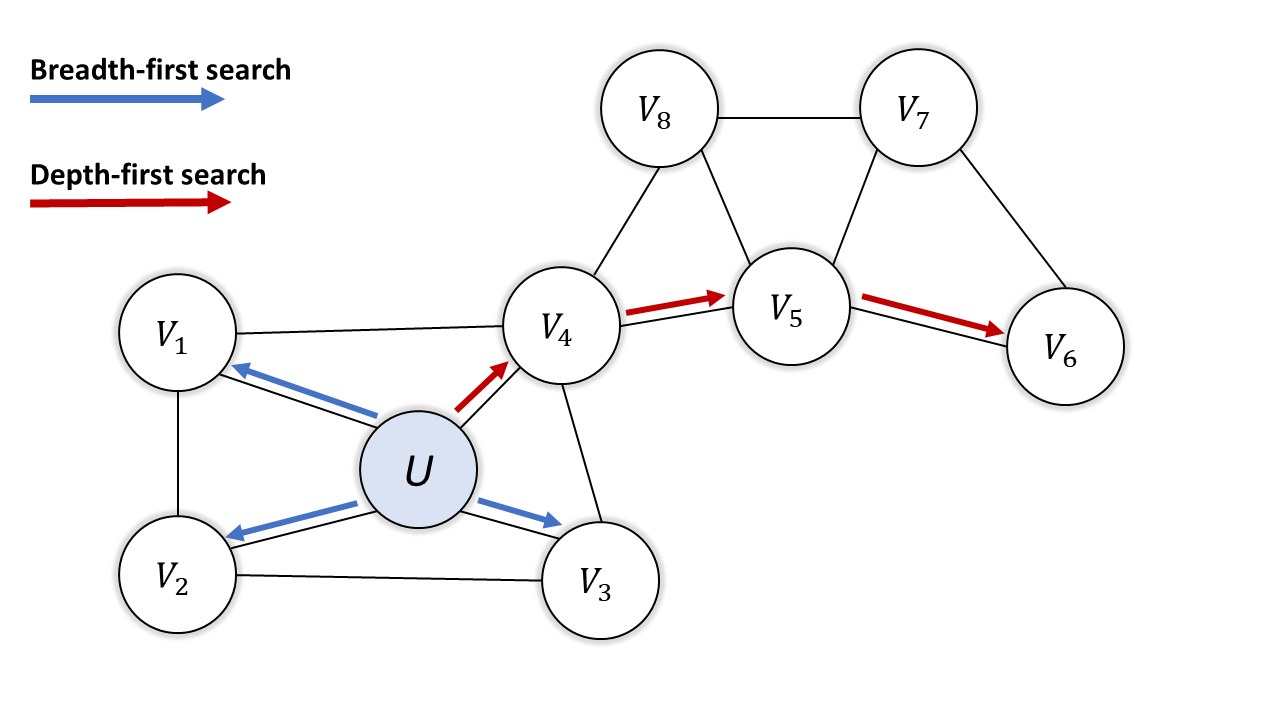}
    \caption{BFS and DFS search strategies from node $u$ ($k = 3$)}
    \label{fig:BFS_DFS}
    \cite{grover2016node2vec}
\end{figure}
\par
\emph{node2vec} introduces two tunable parameters: \textbf{Return parameter $p$} and \textbf{In-Out parameter $q$}. As show in \textbf{Fig. \ref{fig:BFS_DFS}}, these two parameters allow us to move between the two extreme \emph{Breadth-First Search (BFS)} and \emph{Depth-First Search (DFS)} approaches.
\begin{figure}[H]
    \centering
    \includegraphics[height=5.4cm]{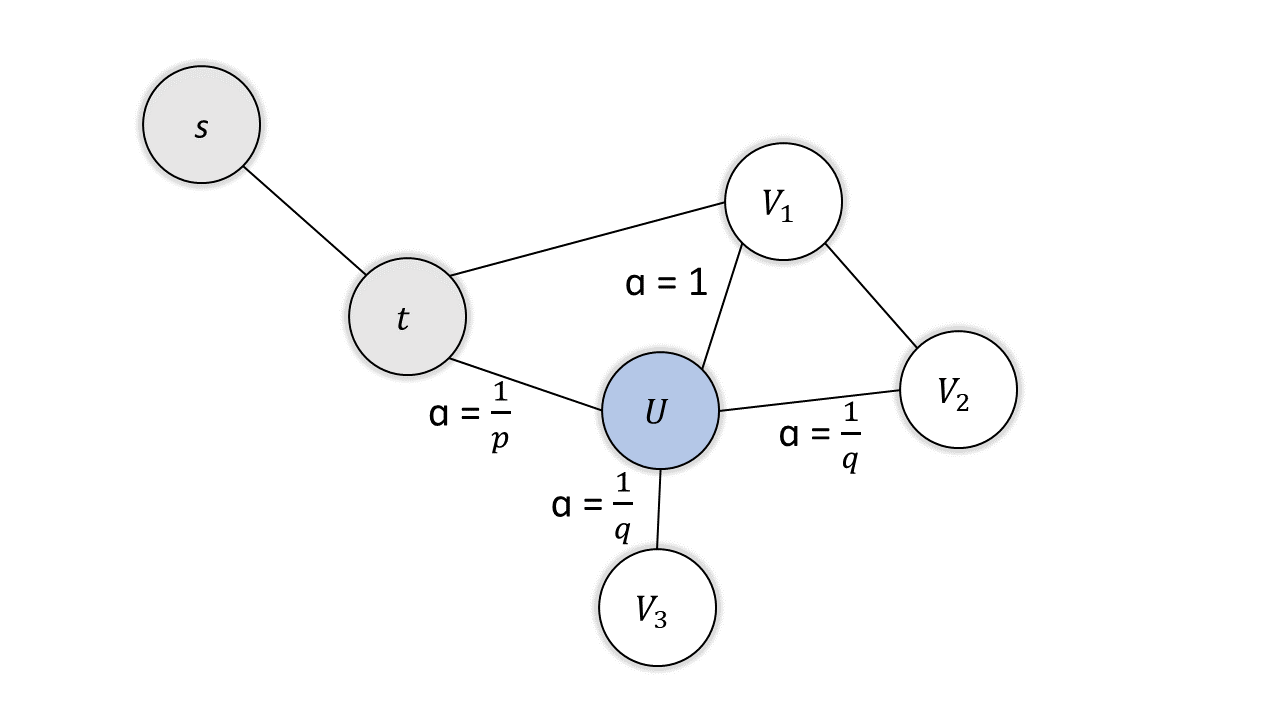}
    \caption{Illustration of the random walk procedure in \emph{node2vec}. The walk just transitioned from $t$ to $v$ and is now evaluating its next step out of node $v$. Edge labels indicate search biases $\alpha$.}
    \label{fig:transition_probability}
    \cite{grover2016node2vec}
\end{figure}
As shown in \textbf{Fig. \ref{fig:transition_probability}}, for a BFS-like walk we need to provide a low value of $p$ and for a DFS-like walk we need a low value of $q$.
\\
\subsubsection{K-Means}
The goal of clustering is to divide the data points into clusters such that the elements assigned to a particular cluster are similar in some predefined sense. Clustering algorithms are divided into two major categories, global and local clustering. In global clustering, every data point is assigned to a cluster in each iteration; whereas in a local approach, the algorithm uses \emph{crawlers} to explore the data points and assign only a subset (usually one) of nodes to a cluster \cite{schaeffer2007graph}.
\par
\emph{K-means} is a  global clustering method originally proposed in 1967 \cite{macqueen1967some}. The main idea is to group points into clusters with the nearest centroids. Algorithms such as K-Means++ \cite{arthur2006k} and Mini-Batch K-Means \cite{sculley2010web} drastically improve K-mean's performance when large-scale data sets are handled. Currently, K-means and methods based on it are widely used for cluster analysis in machine learning applications.
\par
The coarsening scheme introduced in this work is an aggregation based method which uses the preliminary prolongation matrix $\mathbf{\hat{P}}$ (\textbf{Eq. \ref{eq:prolongator}}) as the transfer operator. \emph{node2vec} algorithm is used for the \emph{embedding module}; these embeddings are then utilized by \emph{K-Means} in the \emph{clustering module} to aggregate the nodes.
\par
\emph{K-Means} is a clustering algorithm that tries to divide $N$ samples into $K$ disjoint clusters $C$. K-Means initializes $K$ randomly selected centroids in the sample space and each node is assigned to the closest centroid. The idea is to iteratively update centroid positions so that sum-of-squares of the distances from the corresponding centroids is minimized:
\begin{equation}
    \centering
    argmin_C\sum_{i=1}^{k}\sum_{x \in C_i}||x-\mu_{i}||^2
\end{equation}
\par
\emph{K-Means++} introduces an improvement on initialization algorithm and selecting the initial positions of the centroids. Selecting centroids intelligently increases the speed of \emph{K-Means} clustering algorithm \cite{arthur2006k}.
\par
\emph{K-Means} clustering algorithm does not scale well to large number of samples. Since the graphs that we are working with usually contain more than 10k number of nodes, we are using a customized version of K-Means called \emph{Mini-Batch K-Means} \cite{sculley2010web}.
\par
\emph{Mini-Batch K-Means} applies the same K-Means algorithm but it doesn't run the algorithm on all of the samples at the same instance. It divides the original samples to multiple batches and then runs the algorithm on each batch of samples.
\begin{algorithm}[H]
\caption{Mini-Batch K-Means \cite{sculley2010web}}
\hspace*{\algorithmicindent} \textbf{Input:} $k$, mini-batch size $b$, iterations $t$, data set $X$ \\
\hspace*{\algorithmicindent} \textbf{Output:} Clusters
\begin{algorithmic}[1]
\State $\emph{$v$} \gets 0$
\For{$i = 1$ to $t$}
\State $M\gets b$ samples picked randomly
\For{$x \in M$}
\State $d[x] \gets f(C, x)$ \algorithmiccomment{Cache nearest center to $x$}
\EndFor
\For{$x \in M$}
\State $c \gets d[x]$ \algorithmiccomment{Get cached center for this $x$}
\State $v[c] \gets v[c] + 1$ \algorithmiccomment{Update per-center counts}
\State $\eta \gets \frac{1}{v[c]}$ \algorithmiccomment{Get per-center learning rate}
\State $c \gets (1-\eta)c + \eta x$ \algorithmiccomment{Take gradient step}
\EndFor
\EndFor
\end{algorithmic}
\end{algorithm}
\section{PROPOSED METHOD}
In this paper, we propose \emph{\textbf{GL-Coarsener}} framework. \emph{\textbf{GL-Coarsener}} aims to cluster fine grid based on nodes' neighbors in the underlying graph. These clusters will then be aggregated to form the coarse grid. Then, the residual equation (\textbf{Eq. \ref{eq:residual}}) will be solved on the coarse grid, and finally, the correction in the coarse grid will be transferred to the fine grid and added to the approximate solution.
\\
In the following sections we explain the process step by step.
\subsection{Pre-Processing}
In the AMG context, we are usually dealing with very large sparse matrices. To facilitate working with extensive data and use memory more efficiently, we convert the data to CSR (Compressed Sparse Row) format. In CSR format, only the non-zero elements of the sparse matrix will be stored, and there will be enough information to restore the original matrix elements when needed.
\par
We need to convert matrix $A$ in \textbf{Eq. \ref{eq:residual}} to a graph. Matrix $A$ corresponds to the original grid's underlying graph, with each element in $A$ indicating the weight of edges of the graph. For example, $A_{ij}$ represents the weight of the edge that connects node $i$ to node $j$.
\subsection{Algorithm}
The first step of our proposed AMG algorithm is pre-smoothing. In this step, we will try to reduce high-frequency errors that exist in the approximate solution. To do this, we will use numerical iterative solvers (e.g., Jacobi, Gauss-Seidel). A few iterations as low as 1-5 suffice to reduce high-frequency errors in most cases:
\begin{equation}
    \begin{gathered}
        Av=f
        \\
        Jacobi: v^{(k+1)}=D^{-1}(f-(L+U))v^{(k)}
    \end{gathered}
\end{equation}
where $A=L+D+U$ and $L$, $D$ and $U$ are \emph{lower triangular}, \emph{main diagonal} and \emph{upper triangular} sub-matrices of $A$, respectively.
\par
After performing a few iterations, the iterative solver's convergence rate drops drastically to the point that there is no noticeable improvement in the solution. This significant drop in convergence rate indicates that the iterative solver has reached its limit on the current grid.
\par
That indicates that we need to find a coarser level representation of the system. The system's low-frequency errors will become high frequency at the coarse level, therefore the convergence rate is recovered.
\par
To do this, we will apply graph representation techniques to generate the embedding space of the graph. As said earlier, we will use \emph{node2vec} algorithm. As shown in the \emph{Embedding Module} of \textbf{Fig. \ref{fig:GL_coarsener}}, \emph{node2vec} maps similarities of the graph to the embedding space so that the nodes which are neighbors in the graph, would be closer to each other in the embedding space. As explained in section \textbf{3.2.1}, \emph{node2vec} algorithm generates vectors of size $d$ for each node in order to represent the graph in a \emph{d-dimensional} space.
\par
Afterward, we will feed the newly obtained embedding space to the \emph{K-Means} clustering algorithm. As shown in the \emph{Clustering Module} of \textbf{Fig. \ref{fig:GL_coarsener}}, \emph{K-means} clustering algorithm receives the embedding vectors and the number of clusters as input and will assign similar nodes into the same cluster. In this work, we will use $\frac{n}{5}$ as the number of clusters, meaning that there will be approximately five nodes in each cluster of nodes.
\begin{figure}[hbt!]
    \centering
    \includegraphics[height=7.3cm]{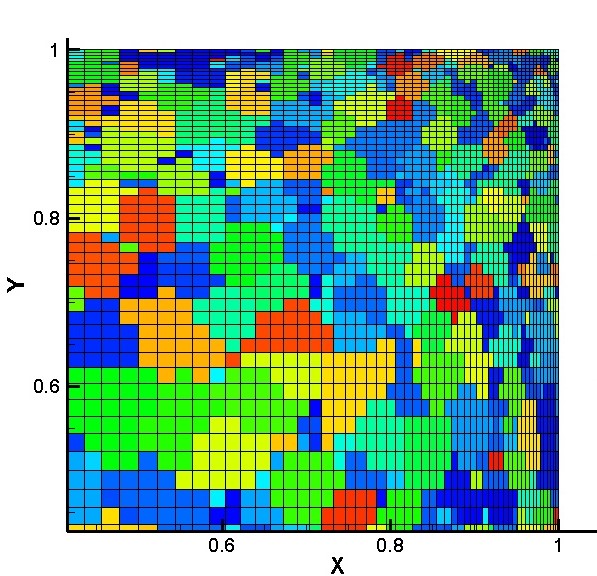}
    \caption{Clusters generated using K-Means clustering method. First, the  SPD system with $10k$ unknowns was generated from a \emph{Poisson} PDE problem and mapped to a $64-dimensional$ embedding space; then using K-Means clustering method, 400 clusters found in the embedding space. The clusters were then colored on the original CFD mesh level. A zoomed level of the top-right corner of underlying grid is shown.}
    \label{fig:cfd_mesh}
\end{figure}
\\
For the example shown in \textbf{Fig. \ref{fig:cfd_mesh}}, \emph{K-Means} clustering algorithm has successfully assigned neighboring nodes of the original grid in to the same clusters.
\par
If our graph size is large (more than 10k nodes), the \emph{K-Means} clustering algorithm might take too long to finish. In this case, as explained above modified version of the algorithm called \emph{Mini-Batch K-Means} is used. As explained in section \textbf{3.2.2}, \emph{Mini-Batch K-Means} accurately approximates the \emph{K-Means} clustering algorithm much faster.
\par
After generating the clusters, we will construct the coarse level coefficient matrix. This step differentiates aggregation-based AMG methods. To do that, we will use \emph{Standard Aggregation} method \cite{vanvek1996algebraic}.
\par
\emph{Standard Aggregation} method proposes a simple approach for constructing the prolongation operator, which will be used to construct the coarse level equation system. The prolongation operator could be either \emph{rough} or \emph{smooth}. A rough prolongation operator could be constructed as shown in \textbf{Eq. \ref{eq:prolongator}}. In small problems, rough aggregation yields a faster convergence.
\newline
Using iterative methods like \emph{Jacobi}, \emph{damped Jacobi}, \emph{Gauss-Seidel} and \emph{SOR}, we could smooth the rough prolongation operator and get a smoother operator:
\begin{equation}
    \begin{gathered}
        P_J = (I - D^{-1}A)\hat{P}
        \\
        P_{\omega J} = (I - \omega D^{-1}A)\hat{P}
        \\
        P_{GS} = (I - (D-L)^{-1}A)\hat{P}
        \\
        P_{SOR} = (I - \omega(D - \omega L)^{-1}A)\hat{P}
    \end{gathered}
\end{equation}
$D$ is the main diagonal of coefficient matrix $A$, $L$ is the strictly lower triangular part of $A$, and $\omega$ is the damping or relaxation parameter.
\par
Once the prolongation operator is constructed, we need to find the residual of system using the approximate solution calculated in pre-smoothing step of the \emph{AMG}. This could be achieved using \textbf{Eq. \ref{eq:residual}}. The residual equation is then transferred to the coarse level:
\begin{equation}
    \begin{gathered}
        A_c = P ^TAP
        \\
        r_c = P^Tr
    \end{gathered}
\end{equation}
\par
This gives the coarse grid equation as:
\begin{equation}
    \begin{gathered}
        A_ce_c=r_c
    \end{gathered}
    \label{eq:coarse_residual}
\end{equation}
\par
Now, we need to solve \textbf{Eq. \ref{eq:coarse_residual}} for $e_c$. $e_c$ on the coarse level corresponds to the correction we need to make on the fine level. Unlike the residual equation in fine level, the system in coarse level features has relatively high-frequency errors. This equation could be solved recursively using the same \emph{AMG} principles.
\par
After solving the residual equation on the coarse level, using the prolongation operator, we need to transfer back the correction calculated on the coarse grid and add the correction to the approximate solution found in pre-smoothing step.
\begin{equation}
    \begin{gathered}
        e = Pe_c
        \\
        u = v + e
    \end{gathered}
\end{equation}
\begin{algorithm}[H]
\caption{One iter. of our proposed \emph{AMG} method}
\begin{algorithmic}[1]
\Function{V-Cycle}{$A_f$, $f_f$, $v_f$}
\State $n \gets \text{length of } f_f$
\State $v_f \gets smoothing(A_f, f_f, v_f)$ \algorithmiccomment{Pre-Smoothing}
\State $\emph{$r_{f}$} \gets \emph{$f_{f}$ - $A_{f}v_{f}$}$
\State $embedding \gets node2vec(A_f)$
\State $clusters \gets MiniBatchKMeans(embedding, \floor*{\frac{n}{5}})$
\State initialize $P$ to be an $n$ by $\floor*{\frac{n}{5}}$ matrix
\For{$i \gets 0$ \textbf{to} $n$}
\For{$j \gets 0$ \textbf{to} $\floor*{\frac{n}{5}}$}
\If{$i \in clusters[j]$}
\State $P_{ij} \gets 1$
\Else
\State $P_{ij} \gets 0$
\EndIf
\EndFor
\EndFor
\State $\emph{R} \gets \emph{P}^T$
\State $A_{c} \gets RAP$
\State $r_{c} \gets Rr_{f}$
\If{\emph{Coarsest Grid Achieved}}
\State $\text{Solve } A_{c}e_{c} = r_{c} \text{ for $e_{c}$}$
\Else
\State $r_c \gets \text{V-CYCLE}(A_c, r_c, e_c)$
\EndIf
\State $e_{f} \gets Pr_{c}$
\State \textbf{return} $v_{f} + e_{f}$
\EndFunction
\end{algorithmic}
\end{algorithm}
\subsection{Code}
The code for our proposed method is available in two forms. \emph{\textbf{GL-Coarsener}} code\footnote{Available at \url{https://github.com/rezanmz/GL-Coarsener}} could be used as an stand-alone module to reduce the size of a large graph. A modular code\footnote{Available at \url{https://github.com/rezanmz/AMG}} for \emph{AMG} is also available that could be used to solve large systems of equations.
\section{NUMERICAL EXPERIMENTS}
This section compares the performance of our proposed \emph{AMG} method with different graph-based and aggregation-based \emph{AMG} methods. These methods include \emph{Beck}'s graph-based method and \emph{Vaněk}'s standard aggregation-based method.
\par
In this section we report implementation of the methods described so far to solve linear systems arising from discretization of the \emph{Poisson}'s equation:
\begin{equation}
    \centering
    \Delta\varphi=f
\end{equation}
where $\Delta$ is the \emph{Poisson}'s operator, $\varphi$ and $f$ are real or complex valued functions.
\par
The discretization is handled by Rayan \cite{sani2010rayan} using Finite Volume Method (FVM) and is treated as the input to the current work. The \emph{Poisson}'s problem is a typical problem occurring in modeling of many physical systems including heat transfer, fluid flow and electrostatics. It is common to solve this problem on grids ranging from a few hundred degrees of freedom (unknowns) to millions and occasionally a few billions (for complex research or industrial problems). 
\subsection{Method Setup}
In what follows the details of settings used in the numerical experiments are described.
\subsubsection{AMG V-Cycle}
The method proposed herein is applicable to all of multi-grid strategies (for example, V, W and F cycles). For demonstration, we will use it in conjunction with V-cycle. Properly setting up the V-Cycle plays a crucial role in its convergence rate. For the pre-smoothing step of the V-Cycle, we use two iterations of \emph{Jacobi} smoothing. The reason for choosing \emph{Jacobi} over other iterative solvers like \emph{Gauss-Seidel} is its natural high potential for parallel application. Our proposed method is also naturally highly parallel and therefore selecting a matching smoother makes the overall set up consistent for massively parallel applications. After each iteration of V-Cycle, we apply seven more \emph{Jacobi} post-iterations; we found it very effective to damp the remaining high-frequency errors.
\par
Since solving large systems directly is computationally very expensive, we will recursively construct coarser grids to the point that we reach a system with 20 unknowns or less. We then solve the system directly.
\par
For a given system of equations, we need to solve the system over multiple iterations of V-Cycle. Evidently, the construction of the prolongation operators happens just on the first iteration of V-cycle and this prolongation operators are reused for the rest of iterations.
\subsubsection{node2vec}
The purpose of \emph{node2vec} embedding module is to map the underlying graph of the system's adjacency matrix to a \emph{d-dimensional} embedding space. $d$ is a hyper-parameter of \emph{node2vec} that needs to be chosen with care. If we choose a too low value for $d$, the embedding space is not large enough to capture all the latent features of the graph and if we choose an excessive number of dimensions for our embedding space, the model takes too long to train and it adversely affects the performance of the method. By running multiple experiments and comparing the results to the  \emph{standard aggregation} method, we found an embedding dimension of 128 to be sufficient to capture latent features of systems with up to 1.5 millions unknowns.
\par
To explore the graph, \emph{node2vec} needs two hyper-parameters $p$ and $q$. In \emph{AMG} aggregation-based methods, we aim to aggregate neighboring nodes in the same cluster, therefore a \emph{local} view of the graph is more desirable than a \emph{global} view.  As described in section \textbf{3.2.1}, a low value $p$ gives us a BFS-like search of the graph resulting in a \emph{local} view. Here we choose \emph{return parameter $p$} to be 0.1 and \emph{in-out parameter $q$} to be 1. However, numerical experiments indicated that the \emph{\textbf{GL-Coarsener}} algorithm is not much sensitive to values of $p$ and $q$. We also choose to run $2 \times AverageDegree$\footnote{In a graph, average degree is simply the average number of in/out edges per node: $Average Degree=\frac{|Edges|}{|Nodes|}$} walks of length 10 to make sure that we capture all the neighbors.
\subsubsection{Mini-Batch K-Means}
To generate the clusters of neighboring nodes in the embedding space, we use \emph{Mini-Batch K-Means}. If we choose the number of clusters to be very low, we will lose too much information when we move between fine and coarse level; and if the number of clusters is very high, in large systems we will have a very long hierarchy of coarser levels; therefore, the runtime increases. By numerical experiments, we found that if we approximately cluster every 5 nodes of the original fine graph into a cluster, we get a good convergence rate; therefore, we choose the number of clusters to be $\frac{n}{5}$. To run the clustering algorithm faster, we choose a batch size of $\frac{n}{15}$.
\subsection{Evaluation}
We evaluate the performance of the \emph{AMG} methods with \emph{V-Cycles} of dynamic depth, meaning that we recursively construct coarser grids to the point we reach a linear system with 20 unknowns or less. The system is then solved directly, and the coarse correction is interpolated back on the respective fine grid. After each iteration, we compute the \emph{infinity} norm of the residual:
\begin{equation}
    \|r\|_\infty := max(|r_1|, ..., |r_n|)
\end{equation}

\subsubsection{Stopping Criteria}
Since \emph{AMG} is an iterative method computing successive approximations to the solution, we should use stopping criteria to determine when to stop iterating. In these experiments, we use:
\begin{equation}
    \|r\|_\infty < 10^{-4}
\end{equation}
as the stopping criterion for the iterative method.
\subsection{Results}
\begin{figure}[H]
    \centering
    \includegraphics[height=6.5cm]{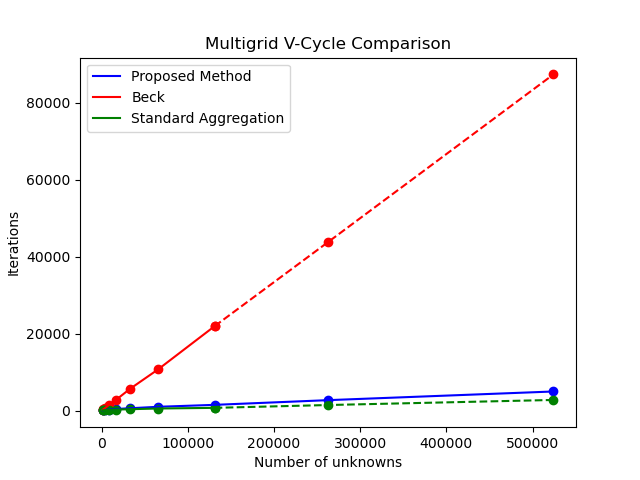}
    \caption{A comparison of different \emph{AMG} methods with stopping criteria of $\|r\|_\infty < 10^{-4}$. For systems larger than 128k unknowns, performances of Beck and Standard Aggregation methods were projected using extrapolation techniques. Graph-based methods (like \emph{Beck}) generally have lower convergence rate than aggregation-based methods (like our methods and \emph{Vaněk}'s standard aggregation)}
    \label{fig:iterations_unknowns}
\end{figure}
\begin{table}[h!]
\centering
\begin{tabular}{||c | c  c  c||} 
 \hline
 System Size & Beck & Vaněk & Proposed Method \\ [0.5ex] 
 \hline\hline
 1k & 186 & 49 & 64 \\ 
 2k & 366 & 80 & 95 \\
 4k & 691 & 110 & 157 \\
 8k & 1419 & 158 & 260 \\
 16k & 2743 & 269 & 412 \\
 32k & 5627 & 356 & 605 \\
 64k & 10709 & 540 & 987 \\ [1ex] 
 \hline
\end{tabular}
\caption{A comparison between multi-grid methods. The data points indicate iterations needed to solve different systems with stopping criteria of $\|r\|_\infty < 10^{-4}$.}
\label{table:1}
\end{table}

As shown in \textbf{Fig. \ref{fig:iterations_unknowns}}, as the size of the system increases, the number of iterations needed for our method increases linearly, closely following \emph{Vaněk}'s standard aggregation-based method and much better than a typical graph based method. Table \ref{table:1} shows the data used to generate \textbf{Fig. \ref{fig:iterations_unknowns}}.
\section{CONCLUSION}
In this paper we propose an aggregation-based \emph{AMG} method. In our proposed method, we use powerful machine-learning and deep-learning techniques to learn the mappings between coarse and fine levels of \emph{AMG}'s V-Cycle. First, we map the system's underlying graph to a \emph{d-dimensional} embedding space, then we use \emph{Mini-Batch K-Means} clustering method to find clusters of neighboring nodes in the system. Using the method provided in \emph{Vaněk}'s standard aggregation,  these clusters are then used to construct the prolongation operator that map the fine level to the coarse level. It has been demonstrated that the GL-Corsener framework while relying on machine learning techniques performs well compared to other existing AMG methods. Another advantage of the method is its naturally parallel characteristic. 
\bibliography{paper.bib}
\bibliographystyle{paper.bst}
\end{document}